\documentclass[12pt]{article}
\usepackage{amssymb}
\usepackage{amsmath}
\usepackage{latexsym}

\def\Z{\mathbb Z}
\def\Q{\mathbb Q}
\def\R{\mathbb R}
\def\C{\mathbb C}

\begin{document}

\newtheorem{theorem}{Theorem}[section]
\newtheorem{lemma}[theorem]{Lemma}
\newtheorem{proposition}[theorem]{Proposition}
\newtheorem{corollary}[theorem]{Corollary}
\newtheorem{remark}[theorem]{Remark}

\title{The asymptotic Schottky problem}
\author{Lizhen Ji\thanks{Partially Supported by NSF grant 
DMS 0604878}, \ 
Enrico Leuzinger}
\date{}
\maketitle

\begin{abstract} 
Let  $\mathcal M_g$ denote the moduli space of compact Riemann surfaces of genus $g$ and
let $\mathcal A_g$ be the space of principally polarized abelian varieties
of (complex) dimension $g$. Let $J:\mathcal M_g\longrightarrow \mathcal A_g$ be the  map which  
associates to a Riemann surface its Jacobian.   The map $J$ is injective,
and the image $J(\mathcal M_g)$ is contained in a proper subvariety of $\mathcal A_g$
when $g\geq 4$.
The classical and long-studied Schottky problem is to 
characterize the
Jacobian  locus $\mathcal J_g:=J(\mathcal M_g)$ in $\mathcal A_g$. 
In this paper we adress a 
large scale version of this  problem posed by  Farb and called the {\em coarse Schottky problem}:
How does $\mathcal J_g$ look ``from far away'', or how ``dense'' is $\mathcal J_g$
in the sense of coarse geometry? 
The coarse geometry of the Siegel 
modular variety $\mathcal A_g$ is encoded in its asymptotic cone  $\textup{Cone}_\infty(\mathcal A_g)$, which
is a Euclidean simplicial cone of (real) dimension $g$. Our main result asserts that the Jacobian
locus $\mathcal J_g$ 
is ``asymptotically large'', or  ``coarsely dense" in $\mathcal A_g$. More precisely,
the subset of $\textup{Cone}_\infty(\mathcal A_g)$ determinded by $\mathcal J_g$ actually coincides with this cone.
The proof also shows that the Jacobian locus of hyperelliptic curves is coarsely
dense in $\mathcal A_g$ as well.
We also study the boundary points of the Jacobian  locus $\mathcal J_g$ in $\mathcal A_g$
and  in the Baily-Borel
and the Borel-Serre compactification. We show that for large genus $g$ the set of boundary points of $\mathcal J_g$ in
 these compactifications is ``small''.
\end{abstract}

\maketitle

\section{Introduction}

The {\it Siegel upper half space} $\mathcal H_g$ is a  Hermitian
 symmetric space of noncompact type which generalizes the Poincar\'e upper half plane:
$$\mathcal H_g:=\{Z\in \mathbb C^{g\times g}\mid Z\  \textup{ symmetric},  \ \textup{Im}\, Z\  \textup{positive definite}\}.
$$
The symplectic group 
$$Sp(g, \mathbb R)=\{x=\left(\begin{array}{ll}
A &  B\\
C & D
\end{array}\right)\in GL(2g, \R)\mid A, B, C, D\in M_n(\R),
x^t J_g x=J_g\}, $$
where $J_g=
\left(\begin{array}{ll}
0 &  I_g\\
-I_g & 0
\end{array}\right)$ and 
$I_g$ is the identity $g\times g$ matrix,
 acts isometrically, holomorphically and transitively on
 $\mathcal H_g$ by
$$
\left(\begin{array}{ll}
A &  B\\
C & D
\end{array}\right)\cdot Z:= (AZ+B)(CZ+D)^{-1}.
$$
As the stabilizer of  the point $iI_g\in\mathcal H_g$ is isomorphic to $U(g)$, one has the identification $\mathcal H_g\cong Sp(g, \mathbb R)/U(g)$.

The {\it  Siegel modular group} $Sp(g, \mathbb Z)$ is an arithmetic subgroup of $Sp(g, \mathbb R)$ which  acts properly discontinuously on $\mathcal H_g$. The corresponding quotient,  
$\mathcal A_g:=Sp(g, \mathbb Z)\backslash \mathcal H_g$ is called the {\it Siegel modular variety},
and can be identified
with the moduli space of  principally polarized abelian varietes (or tori) of complex dimension $g$.
By \cite{Ba1},  $\mathcal A_g$ is a quasi-projective variety
and admits a compactification that is a normal
projective variety. In the following, this compactification is called the
Baily-Borel compactification and denoted by  $\overline{\mathcal A_g}^{BB}$,
in view of the corresponding compactification for general arithmetic Hermitian locally symmetric spaces
constructed in \cite{BB}. 

  The    {\it moduli space $\mathcal M_g$  of compact Riemann surfaces of genus $g>0$}
is a complex (K\"ahler) manifold (or rather orbifold) of dimension $3g-3$. Let 
 $M\in \mathcal M_g$ be  a Riemann surface and consider a symplectic basis $\{A_j,B_j\}$ for the first homology group $H_1(M,\mathbb Z)$ of $M$.  Associated to this basis is a normalized basis $\{\omega_1,\ldots,\omega_g\}$ of holomorphic 1-forms (or abelian differentials
of the first kind)  satisfying $\int_{A_k}\omega_l=\delta_{kl}$. The corresponding {\it period matrix}
$\Pi$ of $M$ is the complex $g\times g$ matrix with entries
$\Pi_{ij}:=\int_{B_i}\omega_j$. Riemann's bilinear relations \cite[p. 232]{GH}
are equivalent to that $\Pi=(\Pi_{ij})$ belongs to
 the Siegel upper half space $\mathcal H_g$. 
 Then $L:=\mathbb Z^g\oplus\Pi\cdot \mathbb Z^g$
is a lattice in $\mathbb C^g$ and the {\it Jacobian} of the Riemann surface 
$M$ is the torus $\mathbb C^g/L$, which turns out to be an abelian variety, i.e.,
it admits the structure of a projective variety.
Moreover, the intersection pairing on homology  $H_1(M, \Z)$ determines a Hermitian bilinear form 
on $\C^g$ with respect to which  the torus $\mathbb C^g/L$
is principally polarized \cite[p. 359]{GH}.

 The choice of a different homology basis of $H_1(M,\mathbb Z)$
  yields a matrix $\Pi'=\gamma\cdot\Pi$ for 
 some $\gamma\in Sp(g, \mathbb Z)$ and hence  a Jacobian in the same isomorphy class of principally polarized abelian varieties. 
 We thus have the well-defined {\it Jacobian (or period) map} 
 $$J:\mathcal M_g\longrightarrow \mathcal A_g$$ which
associates to a Riemann surface $M$  its Jacobian $J(M)$. 
Intrinsically, the Jacobian variety $J(M)$
is equal to $(H^0(M, \Omega^1))^*/H_1(M, \mathbb Z)$,
where $H^0(M, \Omega^1)$ is the space of holomorphic $1$-forms, and the inclusion of $H_1(M, \mathbb Z)$ in the dual space $(H^0(M, \Omega^1))^*$ is obtained
by integrating 1-forms along cycles in $H_1(M, \mathbb Z)$ \cite[p. 36]{GH}. 
By Torelli's Theorem (see \cite[p. 359]{GH}), the Jacobian map $J$ is injective.
The classical \emph{Schottky problem} is to 
characterize the
{\it Jacobian (or period) locus} $\mathcal J_g:=J(\mathcal M_g)$ inside 
the space  $\mathcal A_g$ of all principally polarized abelian varieties.

A lot of work has been done on this important  problem.
Basically  there are two kind of approaches:
(1) the analytic approach, finding equations that ``cut out'' the locus
$J(\mathcal M_g)$ inside $\mathcal A_g$;
(2) the geometric approach, finding geometric properties of a principally
polarized abelian variety that are satisfied only by Jacobians.
For an  nice discussion of the Schottky problem, see
\cite{Mu}.  
More recent surveys of the status of the Schottky problem are \cite{Bea} and \cite{De}.

In \cite{BuS} Buser and Sarnak studied the position of the Jacobian locus $\mathcal J_g$
in $\mathcal A_g$
for large genera $g$.  They consider a certain (systolic) function $m$
which can be thought of as giving a ``distance'' to the boundary of $\mathcal A_g$. 
Then they prove that 
$$\mathcal J_g\subset N_g:=\{x\in\mathcal A_g\mid m(x)\leq\frac{3}{\pi}\log (4g+3)\}.$$
Moreover, 
as $g\to +\infty$, 
$\textup{Vol}(N_g)/\textup{Vol}(\mathcal A_g)=O(g^{-\nu g})$ for any $\nu<1$.
This means that  for large genus the entire Jacobian locus lies in a very small
neigbourhood $N_g$ of the boundary of $\mathcal A_g$.

Motivated by this work of Buser and Sarnak,  B. Farb proposed in \cite[Problem 4.11]{Fa} 
 to study the Schottky problem from the point of view of large scale geometry, called the ``Coarse
Schottky Problem'': How does $\mathcal J_g$ look ``from far away'', or how ``dense" 
 is $\mathcal J_g$ inside $\mathcal A_g$ in the sense of coarse geometry?

This question can be made precise by using the concept of an {\it asymptotic cone} 
(or {\it tangent cone at infinity})
introduced by Gromov.
Recall that a sequence $(X_n,p_n,d_n)$ of unbounded, pointed metric spaces 
converges in the sense of Gromov-Hausdorff, or  {\it Gromov-Hausdorff-converges}, 
to a pointed
metric space $(X,p,d)$ if for every $r>0,  $ the Hausdorff-distance between the balls $B_r(p_n)$ in $(X_n,d_n)$  and  the ball $B_r(p)$ in $(X,d)$ goes to zero as $n\rightarrow\infty$  (see \cite{Gr}, Chapter 3). 
Let $x_0$ be an (arbitrary)  point of  $\mathcal A_g$.
The {\it asymptotic cone} of $\mathcal A_g$ endowed with the locally symmetric metric  $d_{\mathcal A_g}$ 
is defined as the Gromov-Hausdorff-limit of rescaled pointed  spaces:
$$ \textup{Cone}_\infty(\mathcal A_g) := {\mathcal {GH}}-{\lim}_{n\rightarrow\infty} (\mathcal A_g, x_0, \frac{1}{n} d_{\mathcal A_g}).$$
Note that $ \textup{Cone}_\infty(\mathcal A_g)$ is independent of the choice of the base point $x_0$.

  We remark that
in contrast to the case considered here, the   definition of an asymptotic
 cone in general involves the use of ultrafilters, and the limit space may   depend on the
chosen
ultrafilter. Various aspects of  asymptotic cones of general spaces are discussed in Gromov's book \cite{Gr} (see also  \cite{KL}).
In some cases asymptotic cones are  easy to describe. 
For example, the asymptotic cone of the Euclidean space $\R^n$
is isometric to $\R^n$. Similarly, if $C$ is a cone in $\R^n$, then $\textup{Cone}_\infty(C)$
is isometric to $C$.
For another class of examples, 
let  $V$ be a finite volume Riemannian manifold of strictly negative sectional 
curvature and with $k$ cusps, in particular $V$ may be a non-compact, 
finite volume quotient of a rank $1$ symmetric space. Then   $\textup{Cone}_\infty(V)$ is
a  ``cone''  over $k$ points, i.e., $k$ rays with a common origin. 
For  Siegel's modular variety,   $ \textup{Cone}_\infty(\mathcal A_g)$ is known to be isometric
to a $g$-dimensional metric cone over a simplex 
(see Section 2 below).

Farb's question can now be stated as follows \cite[Problem 4.11]{Fa}:

\vspace{1ex}

\noindent{\bf Coarse Schottky problem:}
{\em Describe, as a subset  of a $g$-dimensional metric cone, the subset of
$\textup{Cone}_\infty(\mathcal A_g)$ determined by the Jacobian locus
$\mathcal J_g$ in $\mathcal A_g$}. 

\vspace{1ex}

 Farb also asked to determine the metric distortion of $\mathcal J_g$
inside $\mathcal A_g$ \cite[Problem 4.12]{Fa}. See \S \ref{distortion} below for some comments on that problem. 

 Our first result solves the coarse Schottky problem. It asserts that
the  locus $\mathcal J_g$ is asymptotically  ``dense''. More precisely, we have

\begin{theorem}\label{main} 
Let $\textup{Cone}_\infty(\mathcal A_g)$ 
be the asymptotic
cone of Siegel's modular variety.
Then the  subset of  $\textup{Cone}_\infty(\mathcal A_g)$ determined by the Jacobian locus   $\mathcal J_g\subset \mathcal A_g$  is equal to the entire  $\textup{Cone}_\infty(\mathcal A_g)$.
More specifically, there exists a constant $\delta_g$ depending only on  $g$ such that
$\mathcal A_g$ is contained in a $\delta_g$-neighbourhood of $\mathcal J_g$.
\end{theorem}

In view of the results of Buser and Sarnak Theorem 1.1  might be surprising at first sight.
 Note however that \cite{BuS} deals with the asymptotic
situation when the genus $g\to \infty$, while the genus $g$ is fixed in the present paper. 
The result of Buser and Sarnak implies that the constant $\delta_g\to \infty$. 
A open problem is to find an effective bound on $\delta_g$.

Hyperelliptic curves are special among curves and have been intensively studied
in algebraic geometry.
When the genus $g$ is at least $3$, a generic curve in $\mathcal M_g$
is not hyperelliptic. In fact, denote the   subspace of $\mathcal M_g$ consisting of hyperelliptic curves
by $\mathcal {HE}_g$. 
Then $\dim \mathcal {HE}_g=2g-1$  (see \cite[pp. 255-256]{GH}).
 Since $\dim \mathcal M_g=3g-3$, $\mathcal {HE}_g$ is a proper subvariety of $\mathcal M_g$
 when $g\geq 3$.
Again one can ask about the coarse density
of the image $J(\mathcal {HE}_g)$ in $\mathcal A_g$. The answer is

\begin{theorem}\label{he}
The  subset of  $\textup{Cone}_\infty(\mathcal A_g)$ determined by the 
hyperelliptic Jacobian locus   $J(\mathcal {HE}_g)
\subset \mathcal A_g$  is equal to the entire  $\textup{Cone}_\infty(\mathcal A_g)$. Also,
there  exists a constant $\delta'_g$ depending only on $g$ such that
$\mathcal A_g$ is contained in the $\delta'_g$-neighborhood of $\mathcal {HE}_g$.
\end{theorem}

Siegel's modular variety $\mathcal A_g$ is an arithmetic  Hermitian locally symmetric
space and thus admits several compactifications, which are motivated by various applications (see e.g. \cite{BJ}).
The compactification $\overline{\mathcal A_g}^{BB}$ mentioned above 
is a special case of  the Baily-Borel compactification which exists
for general arithmetic Hermitian locally symmetric spaces (see \cite{BB}). The Baily-Borel compactification is 
a normal projective variety. We denote its
 boundary $\overline{\mathcal A_g}^{BB}-\mathcal A_g$
 by $\partial \overline{\mathcal A_g}^{BB}$.

There is another,  larger compactification of arithmetic locally symmetric
spaces $\Gamma \backslash X$ constructed in \cite{BS},
called the Borel-Serre compactification and denoted by $\overline{\Gamma \backslash X}^{BS}$. It is
 a manifold with corners
and the inclusion $\Gamma \backslash X\hookrightarrow \overline{\Gamma \backslash X}^{BS}$
 is a homotopy equivalence when $\Gamma$ is torsion-free. 
 This Borel-Serre
 compactification has many important applications in topology.
 The basic reason is that $\overline{\Gamma \backslash X}^{BS}$  is a classifying
 space of $\Gamma$ which has the structure of a finite CW-complex and the topology
 of its boundary can be described by the rational Tits building of the associated algebraic group. 
We denote the  Borel-Serre compactification of $\mathcal A_g$ by 
$\overline{\mathcal A_g}^{BS}$ and 
its boundary by
$\partial \overline{\mathcal A_g}^{BS}$.

Since each compactification of $\mathcal A_g$ reflects certain structures or sizes
``near  infinity'', it is natural to consider the boundary points
of the period locus $\mathcal J_g$ in these two compactifications.
Let $\overline{\mathcal J_g}^{BB}$ be the closure of $\mathcal J_g$ in $\overline{\mathcal A_g}^{BB}$,
and $\partial\overline{\mathcal J_g}^{BB}=\overline{\mathcal J_g}^{BB}\cap \partial \overline{\mathcal A_g}^{BB}$
be the boundary of $\overline{\mathcal J_g}^{BB}$ in $\overline{\mathcal A_g}^{BB}$. 
Similarly, let  $\partial\overline{\mathcal J_g}^{BS}$ be the boundary
of $\overline{\mathcal J_g}^{BS}$ in $\overline{\mathcal A_g}^{BS}$.

Our next results show that these boundaries form ``small''' proper subsets when $g$ is large.

\begin{theorem}\label{com1}
When $g=2, 3, 4$, the boundary $\partial\overline{\mathcal J_g}^{BB}$ is equal to the whole
boundary $\partial\overline{\mathcal A_g}^{BB}$ of the Baily-Borel compactification.
For $g\geq 5$, $\partial\overline{\mathcal J_g}^{BB}$ is a proper subvariety of $\partial \overline{\mathcal A_g}^{BB}$.
In fact, it is the union of the Jacobian loci $\mathcal J_k$ of  Riemann surfaces of lower genus $k$,
$k\leq g-1$. 
\end{theorem}

\begin{corollary}\label{com2}
When $g=2, 3$, the boundary $\partial\overline{\mathcal J_g}^{BS}$ is equal to the whole
boundary $\partial\overline{\mathcal A_g}^{BS}$.
For $g\geq 5$,  $\partial\overline{\mathcal J_g}^{BS}$ is a proper subspace of $\partial\overline{\mathcal A_g}^{BS}$
of strictly smaller dimension.
\end{corollary}

Recall that moduli space $\mathcal M_g$ is not compact since there are sequences of compact Riemann surfaces which
degenerate (compare Section 3). The  Deligne-Mumford compactification $\overline{\mathcal M_g}^{DM}$ is a (projective) compactification 
 which is obtained by adding stable Riemann surfaces (see \cite{DM}). In proving Theorem \ref{com1} 
and Corollary \ref{com2} we will use  the following Proposition
\ref{ext} (which may be of independent interest).

\begin{proposition}\label{ext}
 The Jacobian map $J:\mathcal M_g\to \mathcal A_g$ 
extends to a holomorphic and hence algebraic map 
$\overline{J}:\overline{\mathcal M_g}^{DM} \to \overline{\mathcal A_g}^{BB}$.
\end{proposition}

Proposition \ref{ext} also yields  another proof of the result in  \cite{Ba2}
that the closure of $J(\mathcal M_g)$ in $\overline{\mathcal A_g}^{BB}$ is a
subvariety. In fact, that closure is equal to the image $J(\overline{\mathcal M_g}^{DM})$
of a projective variety under an algebraic map. 
 Another consequence of Proposition \ref{ext} is a precise description of the topological 
boundary of $\mathcal J_g$
in $\mathcal A_g$ (see Proposition \ref{interior} below). For $g=2,3$ this in turn
 determines the complement of $\mathcal J_g$ in $\mathcal A_g$
(see Corollary \ref{complement}).

{\small
\vspace{.1in}
\noindent{\bf Acknowledgment.}\ The first author would like to thank Richard Hain for helpful conversations. The second author thanks the
University of Michigan at Ann Arbour for its hospitality and support.}

\section{The coarse geometry of Siegel modular varieties} 

 Coarse fundamental domains, which are usually called fundamental sets,
  for arithmetic groups of semisimple Lie groups
 acting on symmetric spaces of noncompact type are provided by reduction theory (see e.g. \cite{Bo1}). 
 
  In order to describe these fundamental sets in the special case
of $Sp(g, \mathbb Z)$ acting on the Siegel upper half space $\mathcal H_g$,
 we first introduce certain subgroups of 
$Sp(g, \mathbb R)$. We set 
$$
A_g:=\{
\left(\begin{array}{ll}
H &  0\\
0 & H^{-1}
\end{array}\right) \in Sp(g, \R) \mid H\ \textup{positive diagonal}\}\ \ \textup{and}
$$$$
N_g:=\{
\left(\begin{array}{ll}
A &  B\\
0 & ^t\!A^{-1}
\end{array}\right)\in  Sp(g, \R)  \mid A\ \textup{upper triangular, $1$ on diagonal}; A\, ^t\!B=B\, ^t\!A\}.
$$

For $a\in\mathbb R_{> 0}$ we define the {\it Weyl chamber} $ C_a\subset A_g$ as the subset 
of those 
$$
\left(\begin{array}{ll}
H &  0\\
0 & H^{-1}
\end{array}\right)\in 
A_g,\ \ H=\textup{Diag}(h_1,h_2,\ldots,h_g)
$$
which satisfy the inequalities 
$$
\begin{array}{ll}
h_ih_j^{-1}\geq a &\textup{for}\ \  1\leq i<j\leq n\\
h_ih_j \geq a &\textup{for}\ \  1\leq i\leq j\leq n.
\end{array}
$$

Then, for $\omega\subset N_g$  bounded, a    {\it Siegel set in $\mathcal H_g$} is of the form
$$
\mathcal S_{a,\omega}:= \omega C_a\cdot i I_g,
$$
where 
$I_g\in M_g(\C)$ is the identity $g\times g$ matrix and 
$iI_g$ is the chosen base point of $ \mathcal H_g$.
 Note that  $A_g \cdot iI_g$ endowed with the metric induced from $\mathcal H_g$ is a maximal
  totally geodesic
  flat submanifold  of the symmetric space $\mathcal H_g$  and that the Weyl chambers 
  $\mathcal C_{a}=C_{a}\cdot iI_g\subset A_g \cdot i I_g, a>0,$  are Euclidean cones over a simplex. 

The following proposition is a concise version of reduction theory for $Sp(g, \mathbb Z)$; for a proof  see \cite{Si} (or also \cite{Bo1}).

\begin{proposition} \label{red}
There are $a>0$ and $\omega\subset N_g$ as above such that $\mathcal S_{a,\omega}$ is a
fundamental set for $Sp(g, \mathbb Z)$, i.e.

\textup{(1)} $\mathcal H_g=Sp(g, \mathbb Z)\cdot\mathcal S_{a,\omega}$ and

\textup{(2)} the set $\{\gamma\in Sp(g, \mathbb Z)\mid \gamma\cdot \mathcal S_{a,\omega}\cap
\mathcal S_{a,\omega}\neq \emptyset\}$ is finite.
\end{proposition}

We next introduce some additional concepts.
A subset $\mathcal N$ of a metric
space $(X,d)$ is called a ($\delta$--){\it net} if there is a positive  constant
$\delta$ such that  $d(p,\mathcal N)\leq \delta$ for all $p\in X$; in particular the
Hausdorff-distance
between $\mathcal N$ and $X$ is at most $\delta$.
 A map between metric spaces
$f:(X,d_X)\longrightarrow (Y,d_Y)$ is a {\it quasi-isometric embedding} if there are constants $C\geq 1$ and $D\geq 0$ such that
for all $p,q\in X$ one has
$$
C^{-1}d_X(x_1,x_2)-D\leq d_Y(f(x_1),f(x_2))\leq C d_X(x_1,x_2)+D.
$$
If, in addition, the image $f(X)$ is a $\delta$-net in $Y$ for some $\delta>0$, then $f$ is a {\it quasi-isometry}.

Now let $\pi:\mathcal H_g\longrightarrow \mathcal A_g=Sp(g, \mathbb Z)\backslash \mathcal H_g$
 denote the canonical projection and 
 let $\mathcal A_g$ be endowed with the locally symmetric metric such that
$\pi$ is a Riemannian covering. Proposition \ref{red} yields that  $\mathcal A_g=\pi(\mathcal S_{a,\omega})$
and that  $\pi$ is a uniformly bounded finite-to-one map.
 Furthermore  the following metric properties hold.

\begin{proposition} \label{chamber}  Let $\mathcal S_{a,\omega}$ be a Siegel set as in Proposition 
\ref{red}, then
there is $a^*\geq a$ such that $\pi$ restricted to  the Weyl chamber $\mathcal C_{a^*}=C_{a^*}\cdot i I_g\subset 
\mathcal S_{a,\omega}$ is an  isometry. Moreover, $\pi(\mathcal C_{a^*})$ is a net 
in $\mathcal A_g$.
In particular, Siegel's modular variety $\mathcal A_g$ is quasi-isometric
to the Euclidean cone $\mathcal C_{a^*}\cong \pi(\mathcal C_{a^*})$
with the multiplicative constant in the quasi-isometry equal to 1.

\end{proposition}

For a proof (in the general setting of locally symmetric spaces) see \cite[Theorem 4.1, Corollary 4.2]{L1}
and \cite[Lemma 5.10]{JM}.

\vspace{2ex}

The asymptotic cones of   general locally symmetric spaces  of higher rank have been determined by Hattori, Ji-MacPherson and Leuzinger 
(see \cite{Hat}, \cite{JM},  \cite{L1}). Since the asymptotic cone
of a Euclidean cone over a simplex is equal to itself,  we have the following identification for Siegel's modular variety.

\begin{proposition}\label{ascone} Let $\mathcal C_{a^*}$ be as in Proposition \ref{chamber} and let $\pi(\mathcal C_{a^*})$ its isometric image in $ \mathcal A_g$. Then the asymptotic cone, $\textup{Cone}_\infty(\mathcal A_g)$, of $\mathcal A_g$ is isometric to the Euclidean cone $\pi(\mathcal C_{a^*})$,
which is a Euclidean cone over a simplex. 
\end{proposition}

\noindent{\bf Remark}.\ 
The  general result proved in  \cite{Hat}, \cite{JM},  \cite{L1}
is as follows.
Let  $V=\Gamma\backslash G/K$ be a locally symmetric space. Then $\textup{Cone}_\infty(V)$ is isometric to the Euclidean cone over a finite simplicial complex 
$\Gamma\backslash\Delta_{\mathbb Q}$, the
  quotient by $\Gamma$ of the rational Tits building $\Delta_{\mathbb Q}$ of $G$.
In the special case of Siegel modular varieties the quotient of the Tits building is just \emph{one} simplex, and the Euclidean cone over it is isometric to the positive Weyl
chamber $\mathcal C_{a^*}$.  This corresponds to the fact that there is only
one $Sp(g, \Z)$-conjugacy class of minimal $\Q$-parabolic subgroups of $G=Sp(n, \R)$.
Thus 
 $\textup{Cone}_\infty(\mathcal A_g)$ is isometric to a Euclidean cone  $\mathcal C_{a^*}$.

\section{Degenerations of surfaces and period matrices}

It is crucial for our aproach  to obtain information about the image of the Jacobian map $J:\mathcal M_g\longrightarrow \mathcal A_g$ when restricted to certain ``thin parts'' of moduli space $\mathcal A_g$, i.e.,  subsets of $\mathcal M_g$ consisting of Riemann surfaces (endowed with a hyperbolic metric) which contain at least one closed geodesic of length less than some fixed small number (see \cite{L2}
for a precise description of these sets).
The basic philosophy is that while it seems to be difficult  to describe the points in the Jacobian locus
$\mathcal J_g$ completely, we  can describe certain points in its boundary,
and these boundary points in turn allow us to describe the asymptotic cone of $\mathcal J_g$ in $\mathcal A_g$. 
Equivalently, we want to understand the extended Jacobian  map $J: \overline{\mathcal M_g}^{DM}
\to \overline{\mathcal A_g}^{BB}$ from 
 the Deligne-Mumford compactification of $\mathcal M_g$ in Proposition \ref{ext}
 and the intersection of $J(\overline{\mathcal M_g}^{DM})$
 with $\mathcal A_g$.

To this end we discuss in this section the degeneration of Riemann surfaces to singular
surfaces with nodes.
The singular surfaces may be regarded as the union of finitely many compact surfaces
with punctures (the latter identified by the local equation $zw=0$). There are two cases of degeneration in $\mathcal M_g$ depending upon whether the node separates the (singular) surface or not. It is well-known that
these two types of degeneration  yield completely  different  limiting
behaviour in the period locus $\mathcal J_g=J(\mathcal M_g)\subset \mathcal A_g$.

We first  discuss a model for the degeneration of $M$ into two surfaces $M_1$, $M_2$ with genera $g_1,g_2>0$
and joined at a node $p$. 
We choose points $p_1,p_2\in M_1,M_2$ and coordinates $z_i:U_i\rightarrow D$ centered at $p_i$ for $i=1,2$ and $D$ the unit disc in $\mathbb C$. 

Let $S:=\{(z,w,t)\mid zw=t, z,w,t\in D\}$ and let $S_t$ be the fiber for fixed $t$.
Note that when $t=0$, $S_t$ is a singular surface with a nodal point at $(z, w)=(0,0)$,
and when $t\neq 0$, $S_t$ is smooth. 

For $t\in D$ remove the discs $|z_i|<|t|$ from $M_1$ and $M_2$ 
and glue the remaining surfaces by the annulus $S_t$
according to the maps
$$
z_1\mapsto (z_1, \frac{t}{z_1},t),\ \ \ 
z_2\mapsto (\frac{t}{z_2},z_2,t).
$$
This yields an analytic family $\mathcal F\longrightarrow D$ with fibres $M_t, t\neq 0$,
being a compact Riemann surface of genus $g_1+g_2$,
 and $M_0$ a stable Riemann surface (or curve)  with node $p$ (corresponding to $p_1,p_2$), i.e., $M_0$ is a point in the boundary of the Deligne-Mumford compactification of 
 the moduli space $\mathcal M_g$.

We next choose  symplectic homology bases of $M_1$ and $M_2$ to get a symplectic homology basis for $M_t$
in such a way that $\{A_j,B_j\mid 1\leq j\leq g_1\}$ are closed curves in $M_1\cap M_t$ and similarly
$\{A_j,B_j\mid g_1<j\leq g_1+g_2\}$ are closed curves in $M_2\cap M_t$.
For the proofs of the
following proposition and its corollary see \cite[p. 38]{Fay} and \cite[Corollary 3.2]{Fay}.

\begin{proposition} \label{Fay1}
For sufficiently small $t$ there is a normalized basis of abelian differentials
$\omega_1,\ldots,\omega_g$ for $M_t$, holomorphic in $t$, with the following expansions
for $1\leq i\leq g_1$ and $g_1< j\leq g_1+g_2$:
$$
\omega_i(x,t)=\left\{\begin{array}{ll}
 \omega_i^{(1)}(x)+O(t^2)&\ \textup{for}\  x\in M_1\setminus U_1,\\
-t\omega_i^{(1)}(p)\,\omega^{(2)}(x,p)+O(t^2) &\ \textup{for}\ x\in M_2\setminus U_2,
\end{array}\right.
$$
$$
\omega_j(x,t)=\left\{\begin{array}{ll}
 \omega_j^{(2)}(x)+O(t^2)&\ \textup{for}\ x\in M_2\setminus U_2,\\
-t\omega_j^{(2)}(p)\,\omega^{(1)}(x,p)+O(t^2)&\ \textup{for}\ x\in M_1\setminus U_1,
\end{array}
\right.
$$
where the $\omega_i^{(1)}$ form a normalized basis for the abelian differentials on $M_1$, 
the $\omega_j^{(2)}$ form a normalized basis for the abelian differentials on $M_2$
and $\omega^{(1)}(x,y),\omega^{(2)}(x,y)$ are the canonical differentials of the second kind on $M_1,M_2$, which have poles (of order $2$) only along the diagonal $x=y$.
\end{proposition}

\begin{corollary}\label{3.2} The period matrix $\Pi(t)$ of $M_t$ associated to the homology basis
described above satisfies
$$\lim_{t\rightarrow 0}\Pi(t)=\left(\begin{array}{cc}
\Pi_1 &  0\\
0 & \Pi_2
\end{array}\right)
$$
where $\Pi_1$ (resp. $\Pi_2$) is the period matrix of $M_1$ (resp. $M_2$)
with respect to their original homology bases. 
\end{corollary}

\begin{remark}
{\em As explained to the first author by R.Hain, this corollary also follows from
 general results of Schmid on degenerations and 
limits of variations of Hodge structures (see \cite[p. 125]{Ha2}). 
}
\end{remark}

We now turn to the  degeneration of a compact Riemann surface of genus $g$ to a singular 
surface of genus $g-1$ with a single, {\it non-separating} node. The construction is similar
to that of separating nodes, except that we now glue in the annulus via local coordinates
$z_a$ and $z_b$ centered at two distinct points $a,b$ 
on a compact Riemann surface $M$ of genus $g-1$.
Again the resulting surfaces form an analytic family $\mathcal F\longrightarrow D$ with fibres $M_t, t\neq 0$, each being a compact Riemann surface of genus $g$ and $M_0$ a stable Riemann surface.
The node is the identification of $a$ and $b$ in $M_0$ and does not disconnect the 
surface when removed (in contrast to the case considered above).

We choose a symplectic homology basis  $\{A_j,B_j\mid 1\leq j\leq g-1\}$ for $M$ away from the points $a,b$. The surfaces $M_t, t\neq 0$, each have genus $g$ so we need two more loops
$A_g,B_g$ for a homology basis of these surfaces: take $A_g$ as the boundary of the disk
$U_b$, and $B_g$ to run across the handle. One then has the following proposition;
 for the proof see \cite[Corollary 3.8]{Fay}. See also  \cite[Prop. 4.1]{We}.

\begin{proposition} \label{Fay2}
For sufficiently small $t$ the period matrix of $M_t$ has  the following expansion:
$$
\Pi(t)=\left(\begin{array}{cc}
\Pi_{ij}+t\pi_{ij} &  a_i+t\pi_{ig-1}\\
 a_j+t\pi_{g-1j} & -\frac{i}{2\pi}\log t+c_0+c_1t
\end{array}\right)+O(t^2)
$$
where $\Pi=(\Pi_{ij})$ is the period matrix of $M$, $\lim_{t\rightarrow 0}\frac{O(t^2)}{t^2}$
is a finite matrix, and $a_j=\int_a^b\omega_j$.
\end{proposition}

\section{The asymptotic cone of the Jacobian locus}

Before proving the main Theorem \ref{main} we emphasize the following fact.

\begin{lemma} The set $\mathcal D_g$ of diagonal matrices in $\mathcal H_g$ is a
totally geodesic submanifold isometric to a product of $ g$ copies of the  Poincar\'e hyperbolic 
plane $\mathcal H_1$.
\end{lemma}

\emph{Proof}. Consider the map
$$
\Phi: \prod_{k=1}^g Sp(1, \mathbb R)\longrightarrow Sp(g, \mathbb R)$$
given by

$$
\Phi(
\left(\begin{array}{ll}
a_1 &  b_1\\
c_1 & d_1
\end{array}\right), \ldots, 
\left(\begin{array}{ll}
a_g &  b_g\\
c_g & d_g
\end{array}\right)):=
\left(\begin{array}{ll}
D_1 & D_2\\
D_3 & D_4\\
\end{array}\right),
$$
where $D_1,D_2,D_3, D_4\in \mathbb C^{g\times g}$ are {\it diagonal} matrices with entries  $(a_1,\ldots, a_g)$,
$(b_1,\ldots, b_g)$, $(c_1,\ldots, c_g)$, $(d_1,\ldots, d_g)$, respectively.
Clearly, $\Phi$ is an isomorphism of $\prod_{k=1}^g Sp(1, \mathbb R)$ onto its image in 
$Sp(g, \mathbb R)$.
A direct calculation then shows that the orbit of $i I_g\in \mathcal H_g$ under $\Phi(\prod_{k=1}^g 
Sp(1, \mathbb R))$ is the set $\mathcal D_g$ of all diagonal matrices in $\mathcal H_g$  and, moreover, is isometric to the product of $g$ real hyperbolic planes 
$\mathcal H_1\cong Sp(1, \mathbb R)/U(1)\cong SL_2(\mathbb R)/SO(2)$. 
That this embedding is totally geodesic follows for instance 
from the Lie triple criterion (see \cite{He}, IV.7).
\hfill$\Box$

\vspace{1ex }

\subsection{The proof of Theorem \ref{main}}

Given a Riemann surface  $M\in \mathcal M_g$, there are $g-1$  separating curves such the corresponding singular surface
is the union of
$g$  tori   with punctures. Each torus coincides with its own 
Jacobian and corresponds to a point $z_k\in  \mathcal H_1$ (resp. $\mathcal A_1=Sp(1, \mathbb Z)\backslash \mathcal H_1$) for $  k=1,\ldots,g$.
We choose a homology basis
for each torus as above and simultaneously shrink all  $g-1$ separating curves.
 Corollary \ref{3.2}  then implies that there exist period matrices
$\Pi(t)\in \mathcal H_g$ of compact Riemann surfaces in $\mathcal M_g$  such that
$$\lim_{t\rightarrow 0}\Pi(t)=\left(\begin{array}{cccc}
z_1 & 0 & \ldots & 0\\
0 & \ddots &  &\vdots\\
\vdots & &\ddots &0 \\
0 & \ldots & 0 & z_g
\end{array}\right)\in \mathcal D_g\subset \mathcal H_g
$$
for the given  $z_k\in  \mathcal H_1$ for $  k=1,\ldots,g$. Thus every point
in the totally geodesic submanifold  $\mathcal D_g$ (Lemma 4.1) is the limit point of a sequence of period matrices of
surfaces in $\mathcal M_g$  in $\mathcal H_g$.
In order to get the corresponding points in the period locus $\mathcal J_g$ we have to pass
to the quotient $\mathcal A_g=Sp(g, \mathbb Z)\backslash \mathcal H_g$ and in general there will be identifications. However, by Proposition \ref{chamber}, there are no identifications if we restrict to points in the Weyl chamber $\mathcal C_{a^*}\subset \mathcal D_g$
(consisting of certain {\it real}, positive diagonal matrices).
We thus conclude that, given $\delta_1>0$ sufficiently small, the  following holds: for any point $p\in \pi(\mathcal C_{a^*})$ there is a point $p'\in \mathcal J_g$ with $d_{\mathcal A_g}(p,p')<\delta_1$, i.e. 
$$
\pi(\mathcal C_{a^*})\subset \mathcal U_{\delta_1}(\mathcal J_g).
$$
On the other hand, also by Proposition \ref{chamber},
$\pi(\mathcal C_{a^*})\subset \mathcal A_g$ is a net in the modular variety, i.e.
there exists $\delta_2>0$ such that 
$$
\mathcal J_g\subset \mathcal A_g\subset \mathcal U_{\delta_2}(\pi(\mathcal C_{a^*})).
$$
The constant $\delta_g:=\delta_1+\delta_2$ only depends on $g$.
It follows that the Hausdorff distance between  $\mathcal  J_g$ and $\pi(\mathcal C_{a^*})$
is finite:
$$
d_{\mathcal H}(\mathcal  J_g,\pi(\mathcal C_{a^*}))\leq \max\{\delta_1,\delta_2\}.
$$
Consequently, with respect to the recaled metrics $\frac{1}{n}d_{\mathcal A_g}$ the Hausdorff distance 
between $\mathcal  J_g$ and $\pi(\mathcal C_{a^*})$ goes to zero if $n\rightarrow\infty$.
Finally, by Proposition \ref{ascone}, the asymptotic cone of $\mathcal A_g$ 
is isometric to the Euclidean cone $\pi(\mathcal C_{a^*})$. The above estimates thus
imply the claim of Theorem \ref{main} and  the proof is complete.

\subsection{The proof of Theorem \ref{he}}

For each elliptic curve $C$ we fix an origin. Then there exists an involution $\iota$ such that
the origin is a fixed point of $\iota$. There is another fixed point of $\iota$.
The key point is to observe that for $g$ elliptic curves $C_1, \cdots, C_g$
with such fixed involutions,
if we glue them together in a chain along points of involution,
then we get a stable hyperelliptic curve $M_0$.
We can open up these nodes of $M_0$ as in Section 3 while preserving an involution
to get a smooth hyperelliptic curve
$M_t$, i.e., we get a family of curves $M_t$ in $\mathcal {HE}_g$ 
which degenerates to $M_0$. As in the
 proof of Theorem \ref{main} one then shows that the totally
geodesic submanifold $\mathcal D_g$ of $\mathcal H_g$ is contained
in the closure of $J(\mathcal {HE}_g)$,
and the same arguments as above complete the proof of Theorem \ref{he}.

\section{Compactifications of Siegel modular varieties}

Before we prove Theorem \ref{com1} and  Corollary \ref{com2} we briefly  review 
the Baily-Borel and the Borel-Serre compactification of moduli space  $\mathcal A_g$.

\subsection{The Baily-Borel compactification}

First, we describe $\overline{\mathcal A_g}^{BB}$. 
Let $G=Sp(g, \R)$ be the symplectic group with the split $\Q$-structure.
For each  maximal $\Q$-parabolic subgroup $P$ of $G$, there is a Baily-Borel $\Q$-boundary
component of $\mathcal H_g$, denoted by $e^{BB}(P)$, 
which is a Siegel upper half space of lower dimension and constructed as follows. 

Let $P=N_P A_P M_P$ be the Langlands decomposition of $P$  with respect to the maximal
compact subgroup $U(g)$ of $G$. Here  $N_P$ is the unipotent radical of $P$, $A_P$ is the split component, and $M_P$ is a
semisimple Lie group, and furthermore, $A_P$ and $M_P$ are stable under the Cartan involution
associated with $U(g)$.
Then $$X_P=M_P/(U(g)\cap M_P)$$
 is called the {\em boundary symmetric space} associated
with the parabolic subgroup $P$.
It turns out that $X_P$ splits canonically as a product:
$$X_P=X_{h, P}\times X_{\ell, P},$$
where $X_{h,P}$ is a Hermitian symmetric space, and $X_{\ell, P}$ is a homothety section
of a symmetric cone and thus also  called a linear symmetric space. 
The Baily-Borel boundary component associated to the parabolic $P$ is  defined by 
$$e^{BB}(P)=X_{h,P}.$$
The Baily-Borel compactification  $\overline{\mathcal A_g}^{BB}$ is then constructed in two steps:
\begin{enumerate}
\item For every maximal proper $\Q$-parabolic subgroup $P$ of $Sp(g, \R)$, attach the $\Q$-boundary
component $e^{BB}(P)$ to get a partial compactification 
$$\overline{\mathcal H_g}^{BB}=
\mathcal H_g\cup \coprod_{\text{\ maximal $\Q$-parabolic subgroups}\ P} e^{BB}(P).$$
\item Show that $\Gamma=Sp(g, \Z)$ acts continuously on $\overline{\mathcal H_g}^{BB}$
with a compact quotient, which can be given the structure of a  projective variety.
\end{enumerate}

For example,  for every $1\leq k\leq g-1$,  the $\Q$-boundary component $e^{BB}(P_{k,\infty})$
of $\mathcal H_g$ that 
corresponds to the maximal $\Q$-parabolic subgroup
$$P_{k,\infty}=\{\begin{pmatrix} A & 0 & B & n \\
m^t & u & n^t & b \\
C & 0 & D & - m\\
0 & 0& 0 & (u^{-1})^t
\end{pmatrix}\mid
\begin{aligned}
  \begin{pmatrix} A & B\\ C& D\end{pmatrix}\in Sp(k, \R), & u\in GL(g-k, \R) \\
 m, n\in M_{k\times g-k}(\R), &\ \  b\in M_{g-k \times g-k}(\R)
 \end{aligned}
  \}$$
can be identified with the Siegel upper half space $\mathcal H_k$.

The topology of the partial compactification $\overline{\mathcal H_g}^{BB}$
is given  by describing how sequences of interior points converge to boundary points.
The boundary component  $e^{BB}(P_{k,\infty})\cong \mathcal H_k$ is 
attached at the infinity of $\mathcal H_g$ as
$$\{\begin{pmatrix} Z & 0 \\ 0 & +i\infty\end{pmatrix}\mid Z\in 
\mathcal H_k\} 
\cong \mathcal H_k,$$
where a sequence of points $Z_n\in \mathcal H_g$ converges to a point $Z_\infty\in \mathcal H_k$
if and only if when $Z_n$ is written in the form
$Z_n=\begin{pmatrix} Z_n' & Z_n'''\\ (Z_n''')^t & Z_n''\end{pmatrix}$, where $Z_n'\in \mathcal H_k$,
the following conditions are satisfied:
\begin{enumerate}
\item $Z_n'\to Z_\infty$ in $\mathcal H_k$,
\item $\text{Im} Z_n''- (\text{Im} Z_n''')^t (\text{Im} Z_n'')^{-1}(\text{Im}Z_n''')\to +\infty.$ (Note that
for a sequence of real symmetric matrix $y_n\in \C^{n\times n}$, $y_n\to +\infty$ means that for every positive definite
symmetric matrix $A\in \C^{n\times n}$, we have $y_n-A>0$ when $n\gg 1$.)
\end{enumerate}

The action of $Sp(n, \Z)$ on $\mathcal H_g$ extends to a continuous action
on  $\overline{\mathcal H_g}^{BB}$. By the reduction theory for
$Sp(g, \Z)$, it can be shown that every $\Q$-boundary component 
is a translate under $Sp(n, \Z)$ of one of the $e^{BB}(P_{k,\infty})=\mathcal H_k$
described above.

\begin{proposition}\label{BB-A}
The Baily-Borel compactification $\overline{\mathcal A_g}^{BB}$
admits the following disjoint decomposition:
$$\overline{\mathcal A_g}^{BB}=\mathcal A_g\cup \coprod_{k=0}^{g-1} \mathcal A_k,$$
where $\mathcal A_0$ consists of only one point. 
\end{proposition}

{\em Outline of the proof}. 
To prove this result, the crucial point is to observe that though the induced action
of $Sp(n,\Z)$ on $\overline{\mathcal H_g}^{BB}$ is not properly discontinuous,
for each $\Q$-boundary component $e^{BB}(P)$, it  ``effectively" induces a discrete
action on it. Specifically, for the boundary component $e^{BB}(P_{k,\infty})\cong \mathcal H_k$,
two boundary points belong to one orbit of $Sp(g, \Z)$ if and only if
they belong to one orbit of the natural action of $Sp(k, \Z)$ on $\mathcal H_k$.

By the reduction theory for $Sp(g, \Z)$ (compare Proposition \ref{red}),
the parabolic subgroups $P_{k,\infty}$, $k=0, \cdots, g-1$, are representatives
of $Sp(g, \Z)$-conjugacy classes of proper $\Q$-parabolic subgroups of $Sp(n,\R)$.
Combined with the previous paragraph, it implies 

$$Sp(g,\Z)\backslash \overline{\mathcal H_g}^{BB}=
Sp(g,\Z)\backslash \mathcal H_g\cup\coprod_{k=0}^{g-1} Sp(k,\Z)\backslash \mathcal H_k=
\mathcal A_g\cup \coprod_{k=0}^{g-1} \mathcal A_k,$$
which completes the proof. \hfill $\Box$

\vspace{1ex}

For more details of the boundary components and the topology
of $\overline{\mathcal H_g}^{BB}$ and $\overline{\mathcal A_g}^{BB}$ we refer to \cite[pp. 36--37]{Na}.

\subsection{The Borel-Serre compactification}

The Borel-Serre compactification $\overline{\mathcal A_g}^{BS}$ can be constructed as follows.
For {\em every} proper $\Q$-parabolic subgroup $P$ of $Sp(g, \R)$, whether it is maximal or not,
define its boundary component $e^{BS}(P)$ by
$$e^{BS}(P)=N_P\times X_P=N_P\times X_{h, P}\times X_{\ell, P}.$$
We emphasize that for every proper $\Q$-parabolic subgroup $P$ of $G=Sp(g, \R)$,
whether it is maximal or not, 
there is a boundary symmetric space $X_P$, which also splits as the product
$X_P=X_{h, P}\times X_{\ell, P},$ and the Hermitian factor $X_{h,P}$
agrees with the Hermitian factor of a unique
 maximal $\Q$-parabolic subgroup $P_{max}$ containing $P$.

Then $\overline{\mathcal A_g}^{BS}$ is constructed in two steps:

\begin{enumerate}
\item For every  $\Q$-parabolic subgroup $P$ of $Sp(g, \R)$, attach the $\Q$-boundary
component $e^{BS}(P)$ to get a partial compactification 
$$\overline{\mathcal H_g}^{BS}=
\mathcal H_g\cup \coprod_{\text{ $\Q$-parabolic subgroups}\ P} e^{BS}(P).$$
\item Show that $\Gamma=Sp(g, \Z)$ acts continuously and properly on $\overline{\mathcal H_g}^{BS}$
with a compact quotient.
\end{enumerate}

For every maximal $\Q$-parabolic subgroup $P$ of $Sp(g, \R)$, there is clearly a projection
from the Borel-Serre boundary component $e^{BS}(P)$ to the Baily-Borel boundary component
$e^{BB}(P)$. Similarly, for a non-maximal $\Q$-parabolic subgroup $P$,
there is also a projection from $e^{BS}(P)$ to $e^{BB}(P_{\max})$. 
This suggests  the following result (see \cite{BJ} for the proof).

\begin{proposition}\label{compare}
The identity map on $\mathcal A_g$ extends to a continuous, surjective map from
$\overline{\mathcal A_g}^{BS}$ to $\overline{\mathcal A_g}^{BB}$.
\end{proposition}

\section{The boundary points at infinity of $\mathcal J_g$ in compactifications of $\mathcal A_g$}

By Theorem \ref{main}
the Jacobian locus  $\mathcal J_g$ is ``asymptotically  dense'',
 e.g. in the sense that it forms  a net in $\mathcal A_g$.
 In contrast Theorem \ref{com1} asserts that
the boundary of the Jacobian locus $\mathcal J_k$ is ``small'' in the boundary of the Baily-Borel compactification
$\overline{\mathcal A_g}^{BB}$.

In order to prove Theorem \ref{com1}, we  need  extension of the Jacobian map $\overline{J}$
 as stated in  Proposition \ref{ext} 
in the introduction. We then  determine the image of the boundary
$\overline{\mathcal M_g}^{DM}-\mathcal M_g$ in $\overline{\mathcal A_g}^{BB}$
under $\overline{J}$.

\subsection{ The proof of  Proposition \ref{ext}}

The moduli space $\mathcal M_g$ is an orbifold and is covered
by the Teichm\"{u}ller space $\mathcal T_g$, which is a simply connected complex manifold.
The Jacobian map $J:\mathcal M_g \to \mathcal A_g$ can be lifted to a
map $\mathcal T_g\to \mathcal H_g$. Since the boundary of $\overline{\mathcal M_g}^{DM}$
consists of divisors with normal crossing, the Borel extension theorem in \cite{Bo2}
yields the desired extension $\overline{J}$ of the Jacobian map $J$.
\vspace{.1in}

\subsection{ The proof of Theorem \ref{com1}.}

First recall \cite[p. 50]{HaM} \cite{Va} that the
 boundary $\overline{\mathcal M_g}^{DM}-\mathcal M_g$ consists of $[\frac{g}{2}]+1$ divisors,
$D_0, D_1, \cdots, D_{[\frac{g}{2}]}$, where a generic point of $D_0$ is a Riemann surface
of genus $g-1$ with two punctures, and for $k=1, \cdots, [\frac{g}{2}]$,
a generic point of $D_k$ is the union of a Riemann surface of genus $k$ with one puncture
and another Riemann surface of genus $g-k$ with one puncture.

Next note that for each $1\leq k\leq [\frac{g}{2}]$, $\mathcal H_k \times  \mathcal H_{g-k}$ is canonically embedded
into $\mathcal H_g$ by 
$$(Z_1, Z_2)\in \mathcal H_k \times \mathcal H_{g-k} \to \begin{pmatrix} Z_1 & 0\\ 0 & Z_2\end{pmatrix}
\in \mathcal H_g.$$ 
Similarly,  the product $\mathcal A_{g_1} \times \mathcal A_{g_2}$
and hence the product of the Jacobian loci  $\mathcal J_{g_1}\times
\mathcal J_{g_2}$ are also mapped into $\mathcal A_g$ by  finite-to-one
maps. In fact, if $g_1\neq g_2$,  $\mathcal A_{g_1} \times \mathcal A_{g_2}$
is  embedded into $\mathcal A_g$. On the other hand, if $g_1=g_2$,
then the quotient of $\mathcal A_{g_1} \times \mathcal A_{g_2}$ by $\mathbb Z/2$-action
$(x, y)\mapsto (y,x)$ is embedded into $\mathcal A_g$.
More generally, for every proper partition of $g$: $g_1+\cdots+ g_j=g$,
where $g_1, \cdots, g_j\geq 1$, and $j\geq 2$, 
the product $\mathcal H_{g_1} \times \cdots\times \mathcal H_{g_j}$ can  be embedded
into  $\mathcal H_g$, and the product $\mathcal J_{g_1}\times \cdots \times \mathcal J_{g_j}$
is also canonically mapped into $\mathcal A_g$ by a finite-to-one map.
We denote the image of $\mathcal J_{g_1}\times \cdots \times \mathcal J_{g_j}$
in $\mathcal A_g$ by $\mathcal J_{g_1}\times \cdots \times \mathcal J_{g_j}/\sim$.

Recall from \S  3 that for
every generic surface $M_0\in D_k$, $k=0, \cdots, [\frac{g}{2}]$, there is an analytic
family of Riemann surface $M_t$, $t\in D$, such that for $t\neq 0$,
$M_t$ is compact.

Suppose that $k\geq 1$. Then Corollary 3.2 (resp. its obvious generalization) 
 implies that
the limit $\lim_{t\to 0}J(M_t)$ exists and by the above remarks  lies in $\mathcal A_g$.
In particular, the generic points of $\overline{J}(D_k)$  are not contained in the boundary
$\partial \overline{\mathcal A_g}^{BB}$.

Now assume that $k=0$, and let $\overline{M_0}$ be the compact Riemann surface
obtained from $M_0$ by adding the two punctures. Then the genus of $\overline{M_0}$
is equal to $g-1$. Let $\overline{J}(\overline{M_0})$ be the Jacobian of $\overline{M_0}$,
which is a point in $\mathcal A_{g-1}$. Identify $\mathcal A_{g-1}$ with a subset
of the boundary of $\overline{\mathcal A_g}^{BB}$ as in Proposition \ref{BB-A}.
Then $\overline{J}(\overline{M_0})$  canonically determines a boundary point of $\overline{\mathcal A_g}^{BB}$.

By Proposition \ref{Fay2},
the periods of the corresponding analytic family of surfaces $M_t$ are given by
$$\Pi(t)=\left(\begin{array}{cc}
\Pi_{ij}+t\pi_{ij} &  a_i+t\pi_{ig-1}\\
 a_j+t\pi_{g-1j} &- \frac{i}{2\pi}\log t+c_0+c_1t
\end{array}\right)+O(t^2).$$

Note that  $\Pi_{ij}+t\pi_{ij}\in \mathcal A_{g-1}$ converges to a point in $\mathcal J_{g-1}$,
$a_i+t\pi_{ig-1}$ is bounded, 
and Im$(- \frac{i}{2\pi}\log t+c_0+c_1t)\to +\infty$ as $t\to 0$.
By the definition of the convergence of interior points
of $\mathcal A_g$ to the boundary points of  $\overline{\mathcal A_g}^{BB}$ in 
Section 5.1, we see that 
the limit $\lim_{t\to 0} J(M_t)$ exists and is equal to
the boundary point $J(\overline{M_0})$.
This implies that when $\mathcal A_{g-1}$ and hence $\mathcal J_{g-1}$ is identified
with a subset of the boundary $\partial \overline{\mathcal A_g}^{BB}$ as above, then
the boundary $\partial \overline{\mathcal J_g}^{BB}$ of the Jacobian locus contains
$\mathcal J_{g-1}$.

We need to show that the boundary $\partial \overline{\mathcal J_g}^{BB}$ is equal to
the closure of $\mathcal J_{g-1}$ and determine its closure. 
By Proposition \ref{ext}, it is contained in the image of the boundary divisors
$D_0, D_1, \cdots, D_{[\frac{g}{2}]}$ under the extended map $\overline{J}:\overline{\mathcal M_g}^{DM} \to \overline{\mathcal A_g}^{BB}$. The above discussions show that it is 
contained in $\overline{J}(D_0)$. Let $D_0'$ be the open subvariety parametrizing
Riemann surfaces of genus $g-1$ with two punctures.
Then the above discussion shows that $\overline{J}(D_0')$ is equal
to $\mathcal J_{g-1}$ and hence $\partial \overline{\mathcal J_g}^{BB}$ is equal to
the closure of $\mathcal J_{g-1}$. 

If $g\geq 5$, the closure of $\mathcal J_{g-1}$  in $\mathcal A_{g-1}$
is a proper subvariety (see \cite{Ba2} and the discussion in \S 1).
This implies that  the closure of $\mathcal J_{g-1}$
in $\overline{\mathcal J_g}^{BB}$ is a proper subvariety
of $\partial \overline{\mathcal J_g}^{BB}$,
Note that the  arguments in the previous paragraphs and the results in 
Propositions \ref{Fay1} and \ref{Fay2} can be 
generalized to the cases of multiple pinching geodesics
and the closure of $\mathcal A_{g-1}$ in 
the Baily-Borel compactification $\overline{\mathcal A_g}^{BB}$ is the Baily-Borel compactification
of $\mathcal A_{g-1}$.
Then by induction on $g$, we can show that the image $\overline{J}(D_0)$
and hence the boundary $\partial \overline{\mathcal J_g}^{BB}$ is equal to
the union of $\mathcal J_k$, $k\leq g-1$. This completes the proof of the case when $g\geq 5$,
i.e., 
 the second statement of Theorem \ref{com1}.

If $g\leq 4$, then $\mathcal J_{g-1}$ is dense in $\mathcal A_{g-1}$.
Since $\mathcal A_{g-1}$ is dense in $\partial \overline{\mathcal A_g}^{BB}$,
this implies that for $g\leq 4$, $\partial \overline{\mathcal J_g}^{BB}$ is equal to
$\partial \overline{\mathcal A_g}^{BB}$, which proves the first part of Theorem \ref{com1}.

\subsection{ The proof of Corollary \ref{com2}.}

Let $g=2$ or $3$. Then $\mathcal J_g$ is dense in $\mathcal A_g$.
This clearly implies that $\partial \overline{\mathcal J_g}^{BS}$
is equal to $\partial \overline{\mathcal A_g}^{BS}$.
By Proposition \ref{compare}, the boundary $\partial \overline{\mathcal J_g}^{BS}$ is contained
in the inverse image under the map $\overline{\mathcal J_g}^{BS}\to \overline{\mathcal J_g}^{BB}$
of the boundary  $\partial \overline{\mathcal J_g}^{BB}$.

Now let $g\geq 5$.
By Theorem \ref{com1},   $\overline{\mathcal J_g}^{BB}$ is a proper closed
subset of $\partial \overline{\mathcal J_g}^{BB}$.
It follows 
that for $g\geq 5$,   $\partial \overline{\mathcal J_g}^{BS}$
is a proper subset of $\partial \overline{\mathcal A_g}^{BS}$ of strictly smaller dimension.

\begin{remark}
{\em The above arguments miss the case $g=4$.
It is natural to conjecture that in this case,
$\partial \overline{\mathcal J_g}^{BS}$
is also a proper subset of $\partial \overline{\mathcal A_g}^{BS}$ of strictly
smaller dimension.
}
\end{remark}

\section{The interior boundary points of $\mathcal J_g$ in $\mathcal A_g$}

In the previous section we determined the boundary points of the Jacobian locus $\mathcal J_g$ at infinity
of the Siegel modular variety $\mathcal A_g$. In this section we study the closely related 
problem of identifying 
the interior boundary points 
 of $\mathcal J_g$ in $\mathcal A_g$.
More precisely, let $\overline{\mathcal J_g}$ be the closure of $\mathcal J_g$ in $\mathcal A_g$,
and let $\partial \mathcal J_g=\overline{\mathcal J_g}-\mathcal J_g$.

Recall that an abelian variety is  {\em irreducible} if it is not isomorphic to
a product of two abelian varieties of smaller dimensions.
Then the following result on Jacobian varieties is well-known.

\begin{proposition}\label{irreducible} 
For every compact Riemann surface $M\in \mathcal M_g$, its Jacobian 
$J(M)$ is an irreducible principally polarized abelian variety. 
\end{proposition}

We recall the ideas of the proof for convenience. 
By \cite[p. 320]{GH}, every principally polarized 
abelian variety has a Riemann theta-divisor.
As mentioned in the introduction, the Jacobian variety $J(M)$ is  canonically a 
principally polarized abelian variety. Denote  its theta-divisor by $\Theta$.
By a theorem of Riemann \cite[p. 338]{GH}, $\Theta$ is equal to a translate of the image
$W_{g-1}$ of the symmetric power $M^{g-1}$.
Since $M$  and hence $M^{g-1}$ is irreducible,
$\Theta$ is irreducible.

On the other hand, if $J(M)$ is a reducible principally polarized abelian
variety, $J(M)\cong A_1\times A_2$, and $\Theta_1, \Theta_2$ are the Riemann theta-divisors
of $A_1, A_2$, then $\Theta=\Theta_1\times A_2+A_1\times \Theta_2$
and is reducible. This contradiction proves Proposition \ref{irreducible}. 

\vspace{.1in}

Proposition \ref{irreducible} implies that the sets of reducible Jacobians
$\mathcal J_{g_1}\times \cdots \times \mathcal J_{g_j}/\sim$ are disjoint from
the Jacobian locus $\mathcal J_g$.

\begin{proposition}\label{interior}
The closure $\overline{\mathcal J_g}$ of $\mathcal J_g$ in $\mathcal A_g$ admits
the following decomposition:
$$\overline{\mathcal J_g}=\mathcal J_g\cup \coprod_{g_1+\cdots +g_j=g, j\geq 2} 
\mathcal J_{g_1}\times \cdots \times \mathcal J_{g_j}/\sim.$$
Therefore, the interior boundary $\partial \mathcal J_g=\overline{J_g}\setminus J_g$ is given by
$$\partial \mathcal J_g=\coprod_{g_1+\cdots +g_j=g, j\geq 2} 
\mathcal J_{g_1}\times \cdots \times \mathcal J_{g_j}/\sim.$$
\end{proposition} 

{\em Proof.}
It is well-known that a sequence of Riemann surfaces $M_n$  in $\mathcal M_g$  converges
to a stable Riemann surface in the boundary of $\overline{\mathcal M_g}^{DM}$
if and only if a collection of simple, disjoint closed geodesics
on the Riemann surfaces $M_n$ are pinched. 
 Proposition \ref{Fay2} (resp.  a straight forward generalization) 
 and the discussion in the previous section imply that  
if there is a non-separating geodesic,
then the images $J(M_n)$ diverge to  infinity in $\mathcal A_g$.

Therefore, we can assume that all pinching geodesics are  separating.
Then these pinching geodesics determine a partition of $g$: $g_1+\cdots +g_j=g$.
By Corollary 3.2  (resp. its direct generalization) and the discussion in Section 6.2, their corresponding periods $J(M_n)$
converge to a point in the subset $\mathcal J_{g_1}\times \cdots \mathcal J_{g_j}/\sim$
of reducible Jacobians,
and every point in $\mathcal J_{g_1}\times \cdots \mathcal J_{g_j}/\sim$ is the limit of
such a degenerating sequence. 
This proves that $\partial \mathcal J_g$ contains  the union of 
$\mathcal J_{g_1}\times \cdots \times \mathcal J_{g_j}/\sim$
for  every proper partition $g=g_1+\cdots + g_j$. 

In order to show that $\partial \mathcal J_g$ is actually {\em equal} to this union,
we need  the extension result in  Proposition \ref{ext}. 
In terms of the boundary divisors $D_i$ of $\overline{\mathcal M_g}^{DM}$
 and the extended period map $\overline{J}:\overline{\mathcal M_g}^{DM} \to \overline{\mathcal A_g}^{BB}$
 in Proposition \ref{ext}, the boundary  $\partial \mathcal J_g$
is contained in the union $\cup_{k=0}^{[\frac{g}{2}]}\overline{J}(D_k)
\cap \mathcal A_g$. 
As in the proof of Theorem \ref{com1},  let $D_0'$ be the open subvariety parametrizing
Riemann surfaces of genus $g-1$ with two punctures.
Then the arguments there imply that $\overline{J}(D_0')$ is not contained in $\mathcal A_g$.
This implies that the boundary  $\partial \mathcal J_g$
is contained in the union $\cup_{k=1}^{[\frac{g}{2}]}\overline{J}(D_k)
\cap \mathcal A_g$. 

As in the proof of Theorem \ref{com1} again, 
Propositions \ref{Fay1} and  \ref{Fay2} and their generalizations
to the case of multiple pinching curves imply 
that every point in $\cup_{k=1}^{[\frac{g}{2}]}J(D_k)
\cap \mathcal A_g$  is equal to the Jacobian of the stable curve
obtained by pinching only separating simple, disjoint closed curves. 
(The point is that whenever a non-seperating curve of Riemann surfaces is pinched, their
Jacobian varieties will go to the boundary of $\mathcal A_g$.)  
Then by induction on $g$, it can be shown that $\cup_{k=1}^{[\frac{g}{2}]}J(D_k)
\cap \mathcal A_g$ is equal to the union
$\cup_{g_1+\cdots+g_j, j\geq 2}\mathcal J_{g_1}\times \cdots \times\mathcal J_{g_j}/\sim$.
This completes the proof of Proposition \ref{interior}. 

\begin{remark}
{\em
The interior boundary set $\partial \mathcal J_g$ has already been identified in \cite[p. 74]{Mu},
though without a detailed proof. Related results are also hinted in \cite{Ha1}.
We emphasize that without Proposition \ref{ext}, one can only
conclude that $\partial \mathcal J_g$ contains the union $\coprod_{g_1+\cdots +g_j=g, j\geq 2} 
\mathcal J_{g_1}\times \cdots \times \mathcal J_{g_j}/\sim.$
}
\end{remark}

 As pointed out above,
 $\mathcal A_g$ is a Zariski open subset
of the normal projective variety $\overline{\mathcal A_g}^{BB}$.
By \cite{Ba2}, the closure of the image $J(\mathcal M_g)$ in $\overline{\mathcal A_g}^{BB}$
with respect to the regular topology
is an  algebraic subvariety, i.e., the image $J(\mathcal M_g)$
is a quasi-projective variety.  Since $\dim \mathcal M_g=3g-3$ and $\dim \mathcal A_g=g(g+1)/2$,
it follows that for $g=2, 3$ and only for these values of $g$,
 $J(\mathcal M_g)$ is Zariski dense in $\mathcal A_g$.

\begin{corollary}\label{complement} 
When $g=2, 3$, the complement of $\mathcal J_g$ in $\mathcal A_g$ consists
 exactly of the reducible principally polarized abelian varieties. 
\end{corollary}

\noindent{\em Proof.}   When $g=2, 3$, $\mathcal J_g$ is Zariski dense in $\mathcal A_g$.
Then the complement of $\mathcal J_g$ in $\mathcal A_g$ is 
equal to the interior  boundary $\partial \mathcal J_g$, which consists
of reducible Jacobians by Proposition \ref{interior}.
For every proper partition $g_1+\cdots +g_j=g$ with $g=2, 3$ and 
$g_1, \cdots , g_j\leq 2$, it follows that every reducible
 principally polarized abelian variety of dimension $g$ is a reducible
 Jacobian of a stable curve.
 This completes the proof of Corollary \ref{complement}.

\section{Remarks on the distortion of $\mathcal J_g$ inside $\mathcal A_g$}\label{distortion}

Recall that for any path connected subspace $B$ of a geodesic  metric space $(A, d_A)$, there 
is an induced length function $d_B$ on $B$.
Then a function $f:\R_{\geq 0}\to \R_{\geq 0}$ is called a {\em distortion function}
for $B$ inside $A$ if for all pairs of points of $p, q\in B$,
$$ d_B(p, q)\leq f(d_A(p,q)).$$
Clearly, if $B$ is a totally geodesic subspace, then one can take  $f(x)=x$.

In \cite[Problem 4.12]{Fa}, Farb raised the following problem.

\vspace{1ex}

\noindent 
{\em Compute the distortion of the Jacobian locus $\mathcal J_g$ inside $\mathcal A_g$.}

\vspace{1ex}

The expectation in \cite{Fa} is that  any distortion function of $\mathcal J_g$
in $\mathcal A_g$ with respect to the locally symmetric metric is huge and might be exponential.
Based on the results of the previous sections, we think that the distortion of 
$\mathcal J_g$ in $\mathcal A_g$ is of a quite different 
nature.  In fact, we suspect that there are sequences of pairs of points $p_n, q_n\in \mathcal J_g$
such that $d_{\mathcal J_g}(p_n, q_n)$ is bounded away from 0,
or even goes to  infinity, while $d_{\mathcal A_g}(p_n, q_n)$
 goes to 0.
Roughly speaking, when $\mathcal M_g$ is embedded into $\mathcal A_g$, it is folded up,
and for some parts near the boundary of $\mathcal M_g$, 
the different sheets are becoming closer and closer.

To provide such examples of sequences of Riemann surfaces in $\mathcal M_g$, 
we consider the case $g=4$ and start with four distinct Riemann surfaces of genus 1:
$S_1, S_2, S_3, S_4$. Fix them for the moment, though they will move
to  infinity of $\mathcal M_1$ later.
Glue them together to get  analytic families as in \S 3 in the order
$S_1, S_2, S_3, S_4$. In particular, $S_1$ is only connected to $S_2$, but
$S_2$ is glued to both $S_1$ and $S_3$. Similarly, $S_3$ is glued
with both $S_2$ and $S_4$, and $S_4$ is only glued with $S_3$.
Such a family of Riemann surfaces  depends on three parameters $t_1, t_2, t_3\in D$.
We denote these surfaces by $M_{t_1, t_2, t_3}$.

Now we switch the order and glue the surfaces together in the order
$S_2, S_1, S_4, S_3$ so that $S_1$ is now glued with both $S_2$ and $S_4$,
but $S_2$ is only glued with $S_1$. We denote this new family by
$\tilde{M}_{t_1, t_2, t_3}$. 
Due to the different orders,  the compact Riemann surfaces
$M_{t_1, t_2, t_3}$ and $ \tilde{M}_{t_1, t_2, t_3}$ in $\mathcal M_4$
are not isomorphic to each other
for $t_1, t_2, t_3$ sufficiently small.

Now under the extended Jacobian map, the degenerate surfaces
$M_{0,0,0}$ and $\tilde{M}_{0,0,0}$ have the same image:
$$\overline{J}(M_{0,0,0})=\overline{J}(\tilde{M}_{0,0,0}).$$
The crucial point is that these images do not depend on the order of $S_1, S_2, S_3, S_4$.
In fact, they are both equal to the product of the Jacobians $J(S_1), J(S_2), J(S_3), J(S_4)$.
This means that when $t_1, t_2, t_3$  are very small,
$d_{\mathcal A_g}(J(M_{t_1, t_2, t_3}), J(\tilde{M}_{t_1, t_2, t_3}))$ is small.
(Note that they are in the interior of $\mathcal A_g$ and are close to each other.)

On the other hand, in view of Proposition \ref{Fay1}, when $t_1, t_2, t_3\to 0$,
the images $J(M_{t_1, t_2, t_3})$ are basically contained in a Weyl
chamber of $\mathcal A_g$.
Similarly,  for the other family $\tilde{M}_{t_1, t_2, t_3}$,
the images $J(\tilde{M}_{t_1, t_2, t_3})$  are basically contained in another Weyl chamber.
This implies that $d_{\mathcal J_g}(J(M_{t_1, t_2, t_3}), J(\tilde{M}_{t_1, t_2, t_3}))$
is bounded away from 0.
In fact, it is likely that the distance
$d_{\mathcal J_g}(J(M_{t_1, t_2, t_3}), J(\tilde{M}_{t_1, t_2, t_3}))$ goes to infinity.
The reason is that in order to go from one such chamber to another one through
the Jacobian locus $\mathcal J_g$,  there is no shortcut,
and  we need to go through a fixed compact region in $\mathcal A_g$ (or $\mathcal M_g$).
Then Proposition \ref{Fay1} implies the claimed growth of the distance.

The above discussion  indicates that one can  cut $\mathcal M_g$ into finitely many suitable pieces,
whose images in $\mathcal A_g$ under the Jacobian map $J$ have a distortion that is asymptotically
negligible.

\vspace{4ex}

\noindent\textsc{Department of Mathematics\\
 University of Michigan, Ann Arbor, MI 48109 (USA)}

\vspace{1ex}

\noindent{\tt lji@umich.edu}

\vspace{2ex}

\noindent\textsc{Institute of Algebra and Geometry\\
 University of Karlsruhe, 76131 Karlsruhe (Germany)}

\vspace{1ex}

\noindent{\tt enrico.leuzinger@math.uka.de}

\begin{thebibliography}{999}

\bibitem{Ba1}
\textsc{W. Baily}, Satake's compactification of $V_{n}$, 
\emph{Amer. J. Math.} {\bf 80}  (1958), 348--364. 


\bibitem{Ba2}
\textsc{W. Baily}, On the theory of $\theta $-functions, the moduli of abelian varieties, 
and the moduli of curves, \emph{Ann. of Math.}  {\bf 75} (1962) 342--381. 


\bibitem{BB}
\textsc{W. Baily, A. Borel},  Compactification of arithmetic quotients of bounded symmetric domains,
\emph{Ann. of Math.} {\bf 84} (1966),  442--528.


\bibitem{Bea}
\textsc{A. Beauville}, Le probl\`eme de Schottky et la conjecture de Novikov, S\'eminaire Bourbaki
1986/87, \emph{Ast\'erisque} {\bf 152-153} (1988), 101--112.
\bibitem{Bo1}
\textsc{A. Borel}, Introduction aux groupes arithm\'etiques, Paris, 1969.

\bibitem{Bo2}
\textsc{A. Borel}, Some metric properties of arithmetic quotients of symmetric spaces and an extension theorem, 
\emph{J. Differential Geometry} {\bf 6} (1972), 543--560.


\bibitem{BJ}
\textsc{A. Borel, L. Ji}, Compactifications of symmetric and locally symmetric spaces,
\emph{Mathematics: Theory \& Applications}, Birkh\"auser  Boston, 2006.

\bibitem{BS}
\textsc{A. Borel, J.P. Serre}, Corners and arithmetic groups,
 \emph{Comment. Math. Helv.} {\bf 48} (1973), 436--491. 

\bibitem{BuS}
\textsc{P. Buser, P. Sarnak}, On the period matrix of a Riemann surface of large genus
(with an appendix by J.H. Conway and N.J.A. Sloane), \emph{Invent. math.} {\bf 117} (1994), 27--56.  
\bibitem{De}
\textsc{O. Debarre}, The Schottky problem: an update, \emph{Current topics in complex algebraic geometry, Math. Sci. Res. Inst. Publ.} {\bf 28}, Cambridge, 1995, 57--64.

\bibitem{DM}
\textsc{P. Deligne, D. Mumford}, 
The irreducibility of the space of curves of given genus, 
\emph{Inst. Hautes \'Etudes Sci. Publ. Math.} {\bf 36} (1969), 75--109. 

\bibitem{Fay}
\textsc{J.D. Fay}, Theta Functions on Riemann surfaces, \emph{Lecture Notes in Mathematics}
{\bf 352}, Springer, 1973.
\bibitem{Fa}
\textsc{B. Farb},  Some problems on mapping class groups and moduli space, In: Problems on mapping class groups and related topics,  pp. 11--55, 
\emph{Proc. Sympos. Pure Math.} {\bf 74}, Amer. Math. Soc., Providence,  2006.

\bibitem{GH}
\textsc{P. Griffiths, J. Harris}, Principles of algebraic geometry, John Wiley \& Sons, Inc., New York, 
1994. xiv+813 pp.


\bibitem{Gr}
\textsc{M. Gromov}, Metric structures for Riemannian and Non-Riemannian spaces, \emph{Progress in Mathematics}
 {\bf 152}, Birkh\"auser, 1999.
 
 \bibitem{Ha1}
 \textsc{R. Hain},  Locally symmetric families of curves and Jacobians, in {\em
  Moduli of curves and abelian varieties}, 91--108, Aspects Math., E33, Vieweg, Braunschweig, 1999.
 
 
\bibitem{Ha2}
 \textsc{R. Hain},  Periods of limit mixed Hodge structures, in {\em 
  Current developments in mathematics}, 113--133, Int. Press, Somerville, MA, 2003. 
 
 \bibitem{HaM}
 \textsc{J. Harris, I. Morrison}, Moduli of curves. 
 {\em Graduate Texts in Mathematics}, {\bf 187}. Springer-Verlag, New York, 1998. xiv+366 pp.
 
\bibitem{Hat}
\textsc{T. Hattori}, 
 Asymptotic geometry of arithmetic quotients of symmetric spaces,
{\em Math. Z.} {\bf  222 } (1996),  247--277.
\bibitem{He}
\textsc{S. Helgason}, Differential Geometry, Lie Groups, and Symmetric Spaces, \emph{Pure and applied mathematics}
 {\bf 80}, Academic Press, 1978.
\bibitem{JM}
\textsc{L. Ji, R. MacPherson}, Geometry of compactifications of locally symmetric spaces, \emph{Ann. Inst. Fourier, Grenoble} {\bf  52} (2002), 457--559.
\bibitem{KL}
\textsc{B. Kleiner, B. Leeb}, Rigidity of quasi-isometries for symmetric spaces and Euclidean buildings, \emph{Inst. Hautes \'Etudes Sci. Publ. Math.} {\bf  86} (1997), 115-197.
\bibitem{L1}
\textsc{E. Leuzinger},
 Tits Geometry, Arithmetic Groups, and the Proof of a Conjecture of Siegel, \emph{J.  Lie Theory} {\bf  14} (2004), 317--338.
\bibitem{L2}
\textsc{E. Leuzinger},
 Reduction theory for mapping class groups and applications to moduli spaces, preprint 2008.
 
 \bibitem{Mu}
\textsc{D. Mumford}, Curves and their Jacobians, The University of Michigan Press, Ann Arbor, Mich., 
 1975. vi+104 pp. 
 
 \bibitem{Na}
 \textsc{Y. Namikawa}, Toroidal compactification of Siegel spaces,  
\emph{Lecture Notes in Mathematics}, {\bf 812}. Springer, Berlin, 1980. viii+162 pp.
 
 \bibitem{Va}
 \textsc{R. Vakil},  The moduli space of curves and its tautological ring,
 {\em Notices Amer. Math. Soc.}  {\bf 50} (2003), no. 6, 647--658.
 
 
\bibitem{Si}
\textsc{C.L. Siegel}, Collected works, Vol. II, p.108
\bibitem{We}
\textsc{R. Wentworth}, The Asymptotics of the Arakelov-Green's function and Faltings' Delta
invariant, \emph{Commun. Math. Phys.} {\bf 137} (1991), 427--459.
\end{thebibliography}
\end{document}